\documentclass[11pt]{article}
\usepackage{amsmath,amssymb,theorem,amstext,amsgen,amsbsy,amsopn,amsfonts,graphicx,cases}
\usepackage{graphicx}
\usepackage{subfigure}
\usepackage{psfrag}
\usepackage{color}
\usepackage{enumerate}
\usepackage{epstopdf}
\usepackage[numbers,sort&compress]{natbib}
\usepackage{array}

\usepackage{latexsym}
\newtheorem{thm}{Theorem}[section]

\newtheorem{lem}[thm]{Lemma}
\newtheorem{prop}[thm]{Proposition}
\theorembodyfont{\rmfamily}

\def\pf{\bigskip\noindent {\bf Proof.}}

\def\dfn#1{{\sl #1}}

\def\less{\setminus}

\def\pf{\bigskip\noindent {\emph{Proof.}}}

\def\qed{ \hfill\vrule height3pt width6pt depth2pt}

\def\pf{\bigskip\noindent {\bf Proof.  }}

\title{Extremal $H$-free planar graphs}

\author{Yongxin Lan and  Yongtang Shi  \\ 
 Center for Combinatorics and LPMC\\
Nankai University, Tianjin 300071, China\\
and\\
  Zi-Xia Song\thanks{Corresponding author.   \newline E-mail addresses: yxlan0@126.com; shi@nankai.edu.cn; Zixia.Song@ucf.edu}\\
 Department  of Mathematics\\
 University of Central Florida, Orlando, FL 32816, USA
}

\begin{document}
\maketitle
\begin{abstract}
 Given a graph $H$, a graph is \dfn{$H$-free} if it does not contain $H$ as a subgraph. We continue to study the topic of ``extremal" planar graphs, that is,   how many edges
can an $H$-free planar graph on $n$ vertices have?   We define $ex_{_\mathcal{P}}(n,H)$ to be  the maximum number of edges in an $H$-free planar graph on $n $ vertices.     We first obtain several sufficient conditions on $H$ which yield  $ex_{_\mathcal{P}}(n,H)=3n-6$ for all $n\ge |V(H)|$.     We discover that  the chromatic number of $H$ does not play a role, as in   the celebrated  Erd\H{o}s-Stone Theorem.  We then completely determine $ex_{_\mathcal{P}}(n,H)$ when $H$ is a wheel or a star.  Finally, we  
examine the case when $H$ is  a $(t, r)$-fan, that is,   $H$  is isomorphic to  $K_1+tK_{r-1}$, where $t\ge2$ and $r\ge 3$ are integers.   However, determining  $ex_{_\mathcal{P}}(n,H)$,   when $H$ is a  planar subcubic graph,  remains wide open.

\end{abstract}



\baselineskip 16pt
\section{Introduction}
All graphs considered in this paper are finite and simple.  We use $K_t$,  $C_t$ and $P_t$ to denote the complete graph, cycle, and path on $t$ vertices, respectively.  Given a graph $G$, we will use $V(G)$ to denote the vertex set, $E(G)$ the edge set, $|G|$ the number of vertices, $e(G)$ the number of edges, $\delta(G)$ the minimum degree, $\Delta(G)$ the maximum degree. For a vertex $v\in V(G)$, we will use $N_G(v)$ to denote the set of vertices in $G$ which are adjacent to $v$. Let $d_G(v)=|N_G(v)|$ denote the degree of the vertex $v$ in $G$   and $N_{G}[v]=N_{G}(v)\cup \{v\}$. A vertex is a \dfn {$k$-vertex} in $G$ if it has degree $k$. We use   $n_{_k}(G)$ to  denote  the number of $k$-vertices   in   $G$.   For any set  $S \subset V(G)$, 
the subgraph of $G$ induced on  $S$, denoted $G[S]$, is the graph with vertex set $S$ and edge set $\{xy \in E(G) : x, y \in S\}$. We denote by   $G \less S$ the subgraph of $G$ induced  on 
$V(G) \less S$.   If $S=\{v\}$, then we simple write   $G\less v$.  For any two disjoint sets $A$ and $B$ of $V(G)$, we use $e_G(A, B)$ to denote the number of edges between $A$ and $B$ in $G$.  
 The {\dfn{join}} $G+H$ (resp. {\dfn{union}} $G\cup H$)
of two vertex-disjoint graphs
$G$ and $H$ is the graph having vertex set $V(G)\cup V(H)$  and edge set $E(G)
\cup E(H)\cup \{xy\, |\,  x\in V(G),  y\in V(H)\}$  (resp. $E(G)\cup E(H)$). For a positive integer $t$, we use $tG$ to denote disjoint union of $t$ copies of a graph $G$.    Given two isomorphic graphs $G$ and $H$,   we may (with a slight but common abuse of notation) write $G = H$.
  For any positive integer $k$, let  $[k]:=\{1,2, \ldots, k\}$. \medskip 

Given a graph $H$,  a graph is $H$-free if it does not contain   $H$ as a subgraph.   One of the best known results in extremal graph theory is Tur\'an's Theorem~\cite{Turan1941}, which
gives the maximum number of edges that a $K_t$-free graph on $n$ vertices can have. The celebrated  Erd\H{o}s-Stone Theorem~\cite{ErdosStone}  then extends this to the case when $K_t$ is
replaced by an arbitrary graph $H$ with at least one edge, showing that the maximum number of edges possible
is $(1+o(1)){n\choose 2}\left(\frac{\chi(H)-2}{\chi(H)-1}\right)$, where $\chi(H)$ denotes the chromatic number of $H$. This latter result has been called the ``fundamental theorem of extremal graph theory''~\cite{2013Bollobas}.  Tur\'an-type problems when host graphs are hypergraphs are notoriously difficult.  A large quantity of work in this area has been carried out  in determining the maximum number of edges in a $k$-uniform hypergraph  on $n$ vertices without containing  $k$-uniform linear  paths and  cycles (see, for example, \cite{FurediJiang, FurediJiangSeiver, Kostochka1}). Surveys on Tur\'an-type problems of graphs and hypergraphs can be found in \cite{Furedi, Keevash}.\medskip

In this paper, we continue to study the topic of ``extremal" planar graphs, that is,   how many edges
can an $H$-free planar graph on $n$ vertices have?   We define  $ex_{_\mathcal{P}}(n,H)$ to be  the maximum number of edges in an $H$-free planar graph on $n $ vertices.  Dowden~\cite{Dowden2016} initiated the study of $ex_{_\mathcal{P}}(n,H)$ and proved the following results,  where each bound is tight.

 \begin{thm}[\cite{Dowden2016}]\label{Dowden2016}  Let $n$ be a positive integer.
\begin{enumerate}[(a)]
\item $ex_{_\mathcal{P}}(n, C_4)\leq {15}(n-2)/7$ for all $n\geq 4$. \vspace{-8pt}
\item  $ex_{_\mathcal{P}}(n, C_5)\leq (12n-33)/{5}$ for all $n\geq 11$.
\end{enumerate}
\end{thm}

 Let $\Theta_4$ and $\Theta_5$ be obtained from $C_4$ and $C_5$,   respectively,  by  adding  an additional edge joining two non-consecutive vertices.    The present authors~\cite{LSS} studied  $ex_{_\mathcal{P}}(n, H)$ when $H\in\{\Theta_4, \Theta_5, C_6\}$ and when $H$ is a path on at most $9$ vertices. Results from \cite{LSS} are summarized below. 
 \begin{thm}[\cite{LSS}]\label{Theta}  Let $n$ be a positive integer.
\begin{enumerate}[(a)]
 \item $ex_{_\mathcal{P}}(n, \Theta_4)\le {12(n-2)}/5$ for all $n\geq 4$, with  equality when $n\equiv 12 (\rm{mod}\, 20)$. \vspace{-8pt}
 \item  $ex_{_\mathcal{P}}(n, \Theta_5)\le {5(n-2)}/2$ for all $n\ge5$,  with  equality when $n\equiv 50 (\rm{mod}\, 120)$. \vspace{-8pt}
\item $ex_{_\mathcal{P}}(n, C_6)\leq  {18(n-2)}/7$ for all $n\geq6$.   \vspace{-8pt}
\item $ex_{_\mathcal{P}}(n, P_9)\leq   \max\{ {9n}/{4}, {(5n-8)}/{2} \}$. 
 \end{enumerate}
\end{thm}
   
It seems quite non-trivial to determine $ex_{_\mathcal{P}}(n, C_t)$ for all $t\ge7$.  
  In this paper, we first investigate planar graphs $H$  satisfying  $ex_{_\mathcal{P}}(n,H)=3n-6$ for  all  $n\ge|H|$.  This  partially answers a question of   Dowden \cite{Dowden2016}. As observed in \cite{Dowden2016}, for all  $n\ge6$,  
 the planar triangulation $2K_1+C_{n-2}$ is $K_4$-free.  Hence,  $ex_{_\mathcal{P}}(n,H)=3n-6$ for all graphs $H$ which contains $K_4$ as a subgraph and   $n\ge \max\{|H|, 6\}$.  In particular, $ex_{_\mathcal{P}}(n,K_5^-)=3n-6$ for all $n\ge6$, where   $K_p^-$ denotes the   graph obtained from  $K_p$  by deleting one  edge.      Proposition~\ref{main}
   below describes  several sufficient conditions on $K_4$-free planar graphs $H$ such  that $ex_{_\mathcal{P}}(n,H)=3n-6$ for all $n\ge|H|$.

\begin{prop}\label{main}
Let $H$ be a $K_4$-free planar graph and let  $n\ge  |H|$ be an integer. Then $ex_{_\mathcal{P}}(n,H)=3n-6$ if one of the following holds.

\begin{enumerate}[(a)]
\item   $\chi(H)=4$ and $n\ge|H|+2$.\label{g}\vspace{-8pt}
\item $\Delta(H)\ge 7$.\label{a}\vspace{-8pt}

\item $\Delta(H)=6$ and either $n_{_6}(H)+n_{_5}(H)\ge 2$ or  $n_{_6}(H)+n_{_5}(H)=1$ and
      $n_{_4}(H)\ge5$.\label{c}\vspace{-8pt}
\item $\Delta(H)=5$ and  either $H$ has   at least three $5$-vertices or $H$ has exactly two adjacent $5$-vertices. \label{d}\vspace{-8pt}
\item $\Delta(H)=4$ and  $n_{_4}(H)\ge7$. \label{f}\vspace{-8pt}

\item $H$ is $3$-regular  with $|H|\ge9$ or $H$ has at least three vertex-disjoint cycles or $H$ has exactly one vertex $u$ of degree $\Delta(H) \in\{4,5,6\}$  such that   $\Delta(H[N(u)])\ge 3$.\label{e}\vspace{-8pt}
\item   $\delta(H)\ge 4$ or   $H$ has exactly  one  vertex of degree at most $3$.\label{b}

\end{enumerate}

\end{prop}

Proposition~\ref{main} implies that $ex_{_\mathcal{P}}(n,H)=3n-6$ for all $H$ with $n\ge|H|+2$ and either $\chi(H)=4$ or  $\chi(H)=3$ and  $\Delta(H)\ge7$.     Note that $\Theta_4=K_4^-$, and both $K_4^-$  and $K_1+2K_2$ have chromatic number $3$.  Theorem~\ref{Theta}(a)  and   Theorem~\ref{K12K2} (see below)   then demonstrate that  the chromatic number of $H$ does not    play a  role, as in the Erd\H{o}s-Stone Theorem. \medskip

 By Proposition~\ref{main}, $ex_{_\mathcal{P}}(n,H)$ remains unknown for $K_4$-free planar graphs $H$ with   $\Delta(H)=6$, $n_{_6}(H) +n_{_5}(H)=1$ and  $n_{_4}(H)\le4$;
  or   $\Delta(H)=5$ and   $n_{_5}(H)\le2$
   (and the two $5$-vertices are not adjacent when $n_{_5}(H)=2$); or  $\Delta(H)=4$ and  $n_{_4}(H)\le 6$;  or $\Delta(H)\le 3$.  In particular, by Proposition~\ref{main}$(\ref{e})$,  $ex_{_\mathcal{P}}(n,H)$ remains unknown for $K_4$-free planar graphs $H$  with exactly one vertex, say $u$, of degree $\Delta(H)\le6$ and $\Delta(H[N(u)])\le2$.  It seems non-trivial  to determine $ex_{_\mathcal{P}}(n,H)$ for all such $H$.  We next study $ex_{_\mathcal{P}}(n, W_k)$, where
    $ W_k:=K_1+C_{k}$  is a \dfn{wheel} on $k+1\ge 5$ vertices.
 Unlike the classic Tur\'an number of $W_k$ (see \cite{Dzido2013, DJ2017} for more information), the planar Tur\'an number  of $W_k$ can be completely determined. We establish  this  in Theorem~\ref{wheel}.

\begin{thm}\label{wheel}
Let $n, k$  be integers with $n \ge k+1\ge 5$. Then
\[ex_{_\mathcal{P}}(n,W_k)=
\begin{cases}
\, 3n-6 & \text{ if  } \,  k\ge 6,  \text{ or }  k=5 \text{  and }  n\ne 7,  \text{  or } k=4 \text{ and } n\ge12 \\[2mm]

\, 3n-7 & \text{ if  }  \, k=4 \text{  and } \ n\in\{5, 6\}, \text{  or }  k=5\text{  and } n=7\\[2mm]

\, 3n-8 & \text{ if  }  \,   k=4    \text{ and } 7\le n\le 11.
\end{cases}\]
\end{thm}

 A graph is a \dfn{$(t, r)$-fan} if it is isomorphic to $K_1+tK_{r-1}$, where $t\ge2$ and $r\ge 2$ are integers. The classical Tur\'an number   of  $(t, r)$-fan, namely, $ex(n, K_1+tK_{r-1})$,   has also been studied  when $n$ is sufficiently large (see  \cite{EFGG1995, CGPW2003}  for more information).     We next study $ex_{_\mathcal{P}}(n,H)$  
 when $H$ is a $K_4$-free     $(t, r)$-fan, in particular,   when $H\in\{K_1+2K_2, K_{1, t}, K_1+3K_2\}$.   Theorem~\ref{K12K2} below establishes a sharp  upper bound for  $ex_{_\mathcal{P}}(n,K_1+2K_2)$, and Theorem~\ref{star} completely determines the  value of $ex_{_\mathcal{P}}(n,K_{1, t})$ for all $t\ge3$.  However, the upper bound for $ex_{_\mathcal{P}}(n,K_1+3K_2)$ when $n\ge 15$ in Theorem~\ref{K13K2} is not tight.

\begin{thm}\label{K12K2} Let $n\ge5$ be an  integer. Then
 $$2n-3\le ex_{_\mathcal{P}}(n,K_1+2K_2)\le \frac{19n}{8}-4.$$  Furthermore, $ex_{_\mathcal{P}}(n,K_1+2K_2)= \frac{19n}{8}-4$ if and only if  $n$ is divisible by $8$.
 \end{thm}

\begin{thm}\label{star} Let $n, t$ be    integers with $n\ge t+1\ge4$. Then
\[ex_{_\mathcal{P}}(n,K_{1, t})=
\begin{cases}
\, 3n-6 & \text{ if }  \,  t\ge7,     \text{ or } \,  t=6 \text{ and }   n\in\{7,8,9, 10,12\} \\[2mm]
\, 3n-7 & \text{ if } \,  t=6  \text{ and }\,   n=11 \\[2mm]
\, 3n-8 & \text{ if }   \,  t=6   \text{ and } \,   n\in\{13,14\},    \text{ or }\,  t=5  \text{ and }  \,  n=7 \\[2mm]
\, \left\lfloor\frac{(t-1)n}{2}\right\rfloor   & \text{ if  }\,  t\in\{3,4\},  \text{ or }\,  t= 5 \text{ and } \, n\ne 7,   \text{ or } \, t=6 \text{ and } \, n\ge 15.
\end{cases}\]

\end{thm}

\begin{thm}\label{K13K2} Let $n\ge7$ be an  integer. Then $$\left\lfloor\frac{5n}{2}\right\rfloor \le ex_{_\mathcal{P}}(n,K_1+3K_2)< \frac{17n}{6}-4 $$ for all $n\ge 15$ and
\[ex_{_\mathcal{P}}(n,K_1+3K_2)=
\begin{cases}
\, 3n-6 & \text{  if }\,     n\in\{7,8,9, 10,12\} \\[2mm]
\, 3n-7 & \text{  if } \,    n=11 \\[2mm]
\, 3n-8 & \text{  if }\,     n\in\{13,14\}.
\end{cases}\]

\end{thm}
\bigskip

  We need to introduce more notation.  
  Given a plane graph $G$ and an integer $i\ge 3$, an $i$-face in $G$ is a face of order $i$.  Let $f_i$ denote the number of $i$-faces in $G$. Let $T_n$ denote a plane triangulation on $n\ge3$ vertices, and let $T_n^-$ be obtained from $T_n$  with one edge removed. 
   For every integer $n\ge 5$, let $O_n $ 
   denote the unique  outerplane graph with $2n-3$ edges,    maximum degree $4$, and the  outer face of order $n$;   let  $O_n'$ be a different drawing of $O_n$ with  one unique inner face of order $n$; and 
 finally,   let $O^*_n$   be the planar triangulation obtained from $O_n$ and $O'_n$ by identifying the outer face of $O_n$ with the unique $n$-face    of $O_n'$ in such a way that  $ O^*_n$ is a simple graph.  
 The proof of Proposition~\ref{main} relies heavily on the fact that    $O^*_n$, $K_1+O_{n-1}$ and $2K_1+C_{n-2}$ are planar triangulations.  
\medskip

Finally, we shall make use of the following lemma  in the proof of Theorem~\ref{star} and Theorem~\ref{K13K2}.

\begin{lem}[\cite{SH1977}]\label{Grunbaum1}
There does not exist a $4$-regular planar graph on $7$ vertices, or a $5$-regular planar graph on $14$ vertices, or a planar graph on $n\in\{11, 13\}$ vertices with exactly one vertex of degree $4$ and $n-1$ vertices of degree $5$.
\end{lem}

\section{Proof of Proposition~\ref{main}}

Let $H$ and $n$  be given as in the statement.  To prove $(\ref{g})$, assume $\chi(H)=4$ and $n\ge |H|+2$. Since the planar triangulation $2K_1+C_{n-2}$ has no subgraph on $|H|$ vertices with chromatic number $4$, we see that  $2K_1+C_{n-2}$ is $H$-free. Hence,    $ex_{_\mathcal{P}}(n,H)=3n-6$ when  $\chi(H)=4$ and $n\ge |H|+2$.\medskip

To prove $(\ref{a})$,  assume  that  $\Delta(H)\ge7$.
Then  $n\ge |H|\ge8$, and    the   planar
triangulation $O^*_n$ is $H$-free because $\Delta(O^*_n)=6$. Hence, $ex_{_\mathcal{P}}(n,H)=3n-6$  for all $n\ge |H|$.  \medskip

To prove $(\ref{c})$,  assume  $\Delta(H)=6$. Then $n\ge |H|\ge7$.   Assume  first $n_{_6}(H)+n_{_5}(H)\ge 2$.  
Let  $x,y\in V(H)$ be such that $d_H(x)=6$ and $d_H(y)\ge5$.  Then  the   planar
triangulation $K_1+O_{n-1}$ is  $H$-free when $xy\notin E(H)$,  and  the planar triangulation $2K_1+C_{n-2}$ is $H$-free when $xy\in E(H)$. Hence,
 $ex_{_\mathcal{P}}(n,H)=3n-6$ when $n_{_6}(H)+n_{_5}(H)\ge 2$.  Next assume that $n_{_6}(H)+n_{_5}(H)=1$  and $n_{_4}(H)\ge5$.  Then  $n_{_5}(H)=0$. Let $S\subset V(H)$ be the set of all vertices $v$ such that   $d_H(v)=6$ or $d_H(v)=4$. Then  the planar triangulation $2K_1+C_{n-2}$ is $H$-free when $n_{_4}(H)\ge6$ or     $H[S]\neq 2K_1+P_4$,  and  the planar triangulation $ O^*_{n}$ is $H$-free when $H[S]=2K_1+P_4$. It follows that 
 $ex_{_\mathcal{P}}(n,H)=3n-6$ for all $n\ge |H|$.  This proves $(\ref{c})$.  \medskip

To prove $(\ref{d})$, assume $\Delta(H)=5$ and $n_{_5}(H)\ge2$. Then $n\ge |H|\ge 6$.   Let $u, v$ be two distinct $5$-vertices in $H$.  Then either   $n_{_5}(H)\ge 3$ or $n_{_5}(H)=2$ with $uv\in E(H)$.    Note that the planar triangulation  $2K_1+C_{n-2}$ has exactly two non-adjacent  vertices of degree at least $5$   when $n\ge 7$ and  has maximum degree $4$   when $n=6$. Hence,  $2K_1+C_{n-2}$ is $H$-free, and so   $ex_{_\mathcal{P}}(n,H)=3n-6$ for all $n\ge |H|$.
This proves $(\ref{d})$.\medskip

To prove  $(\ref{f})$ and $(\ref{e})$.  Assume $\Delta(H)=4$ and  $n_{_4}(H)\ge7$, or  $H$ is $3$-regular  with $|H|\ge9$,  or $H$ has at least three vertex-disjoint cycles, or $H$ has exactly one vertex $u$ of degree $\Delta(H) \in\{4,5,6\}$  such that   $\Delta(H[N(u)])\ge 3$.  Then the planar triangulation  $2K_1+C_{n-2}$ is $H$-free. Hence,   $ex_{_\mathcal{P}}(n,H)=3n-6$ for all $n\ge |H|$.  \medskip

It remains to prove  $(\ref{b})$.   Assume $\delta(H)\ge4$ or $H$ has exactly  one vertex of degree at most $3$.  Then $n\ge|H|\ge 5$.  Note that   the   planar
triangulation $K_1+O_{n-1}$ is $3$-degenerate and every subgraph of $K_1+O_{n-1}$ has at least two vertices of degree at most $3$,  because every subgraph of $O_{n-1}$ has at least two vertices of degree at most $2$.   Hence, $K_1+O_{n-1}$ is $H$-free,  and so  $ex_{_\mathcal{P}}(n,H)=3n-6$ for all $n\ge |H|$.  This completes the proof of Proposition~\ref{main}. \qed\medskip

\section{Proof of Theorem~\ref{wheel}}

\begin{figure}[htbp]
\centering
\includegraphics*[scale=0.3]{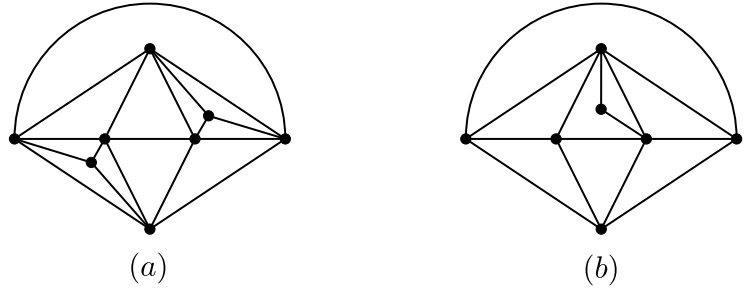}
\caption{$W_k$-free plane graphs with  $n=k+2$ vertices and $3n-12+k$ edges, where   $k\in\{5, 6\}$. }\label{T8T7-}
\end{figure}

\begin{figure}[htbp]
\centering
\includegraphics*[scale=0.32]{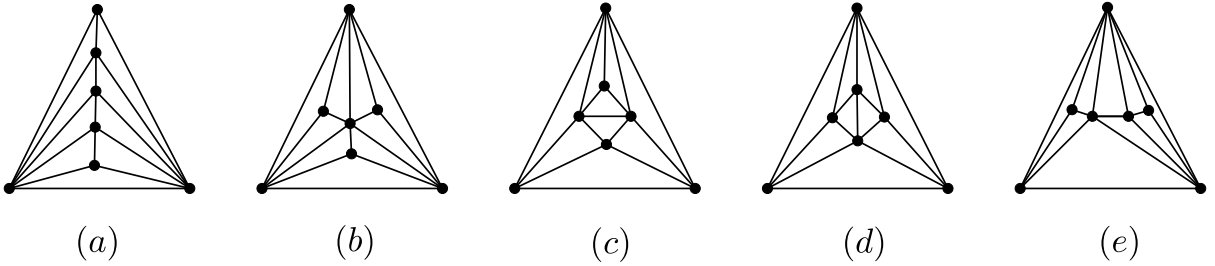}
\caption{All plane triangulations on $7$ vertices. Furthermore, each $T_7$ has a copy of $W_4$ and $W_5$, and each $T_7^-$ has a copy of $W_4$.}\label{seven}
\end{figure}

Let $n, k$ be given as in the statement.
Assume $k\ge7$. By Proposition~\ref{main}$(\ref{a})$, $ex_{_\mathcal{P}}(n,W_k)=3n-6$ for all $n\ge k+1$. Assume next $k\in\{5,6\}$. Since the planar triangulation $2K_1+C_{n-2}$ is $W_k$-free when   $n\ge k+3 $ or $n=k+1$, we see that
$ex_{_\mathcal{P}}(n,W_k)=3n-6$ when    $n\ge k+3 $ or $n=k+1$. We next  determine $ex_{_\mathcal{P}}(n,W_k) $ when $n=k+2$.
For  $k=6$  and $n=8$, the plane triangulation  on $8$ vertices depicted in Figure~\ref{T8T7-}(a) is $W_6$-free and so   $ex_{_\mathcal{P}}(n,W_6)=3n-6$ when $n=8$.
For  $k=5$  and $n=7$,    note that the plane graph with  $7$ vertices and  $14$ edges  given in Figure~\ref{T8T7-}(b) is $W_5$-free. Thus,  $ex_{_\mathcal{P}}(7,W_5)\ge 3\cdot 7-6-1=14$. On the other hand,  all plane triangulations on  $7$  vertices are depicted in Figure \ref{seven}, each  containing a copy of $W_5$.  Hence,  $ex_{_\mathcal{P}}(n,W_5)=3n-7$ when $n=7$.

\begin{figure}[htbp]
\centering
\includegraphics*[scale=0.26]{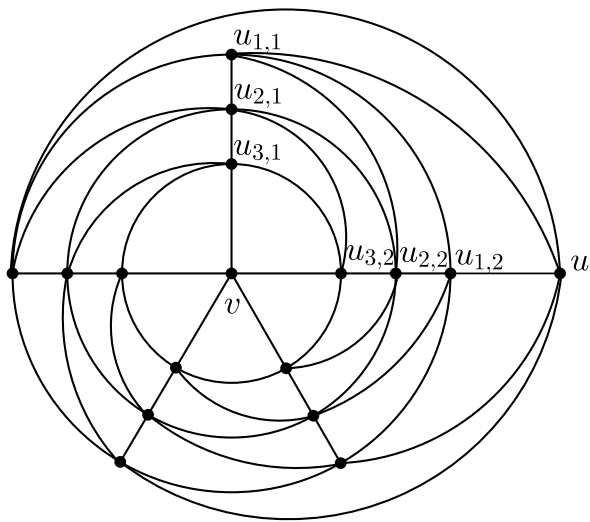}
\caption{Plane triangulation $L_{t}$ when   $t=3$.}\label{Tkt}
\end{figure}

It remains to consider the case when   $k=4$. To show $ex_{_\mathcal{P}}(n,W_4)=3n-6$ for all $n\ge 12$, assume first that    $n=5t+2$ for some  integer $t\ge 2$.  Let $L_t$ be a plane triangulation on $n=5t+2$ vertices constructed as follows:  for each $i\in[t]$, let $C^i $ be a cycle with vertices $u_{i,  1}, u_{i,  2}, \ldots,  u_{i,  5}$ in order. Let $L_t $ be the plane triangulation obtained from disjoint union of $ C^1, \ldots,  C^t $ by adding edges $u_{i,  j}u_{i+1,  j}$ and $u_{i,  j}u_{i+1,  j+1}$   for all $ i\in [t-1]$ and $j\in[5]$, where all arithmetic on the   index $j+1$ here   is done modulo  $5$,  and finally adding two new non-adjacent vertices $u$ and $v$ such that $u$ is  adjacent to all vertices of $C^1 $ and $v$ is adjacent to all vertices of $C^t$. The graph $L_t $ when   $t=3$ is depicted in Figure~\ref{Tkt}. It is worth noting  that $L_t$ is $K_4$-free,  $d_{{L_t}}(u)=d_{{L_t}}(v)=5$, $d_{{L_t}}(u_{1,\, j})=d_{{L_t}}(u_{t,\, j})=5$, $d_{{L_t}}(u_{i,\, j})=6$ for $2\le i\le t-1$ and $j\in[5]$. Furthermore, the subgraph induced by the neighborhood of each vertex in $L_t$ is isomorphic to either $C_5$ or $C_6$. Hence, $L_t$ is $W_4$-free and so $ex_{_\mathcal{P}}(n,W_4)= 3n-6$ when $n=5t+2$ for some  integer $t\ge 2$.\medskip

Next assume  that $n= 5t+2+i$ for some  $i\in[4]$, where   $t\ge2$ is an integer. Note that the plane triangulation $L_t$ on $5t+2$ vertices constructed above contains at least four  pairwise  vertex-disjoint faces.  Let $ F_1, \ldots, F_i $ be any $i$  pairwise  vertex-disjoint faces of $L_t$,  and let $L_t^i$ be the plane triangulation obtained from $L_t$ by adding one vertex, say $x_j$, of degree $3$ to each   $ F_j$ for all $j\in[i]$.
Clearly,  $L_t^i$ is a plane triangulation on $n=5t+2+i $ vertices. By the choice of $ F_1, \ldots, F_i $, we see that $x_1, \ldots, x_i$    are pairwise non-adjacent  in $L_t^i$,  and no two of $x_1, \ldots, x_i$ have common neighbors in  $L_t^i$.
We next show that    $L_t^i$ is $W_4$-free for all $i\in[4]$. Suppose that $L_t^i$ contains a copy of $W_4$ for some $i\in[4]$. Let $H$ be a $W_4$ in $L_t^i$. Then $H$  must contain exactly one, say $x_1$,  of $x_1, \ldots, x_i$,  because $L_t$ is $W_4$-free, and no two  of $x_1, \ldots, x_i$  are adjacent  or have common neighbors in  $L_t^i$.  Let $y, z\in V(H)$ be the two neighbors of $x_1$ such that $yz\notin E(H)$. By the choice of $x_1$, we see that $yz\in E(L_t)$.  But then $L_t[V(H\less x_1)]=K_4$ and so
     $L_t$  contains $K_4$ as a subgraph, a contradiction.  Therefore, $ex_{_\mathcal{P}}(n,W_4)=3n-6$ for all $n\ge 12$.

\begin{figure}[htbp]
\centering
\includegraphics*[scale=0.3]{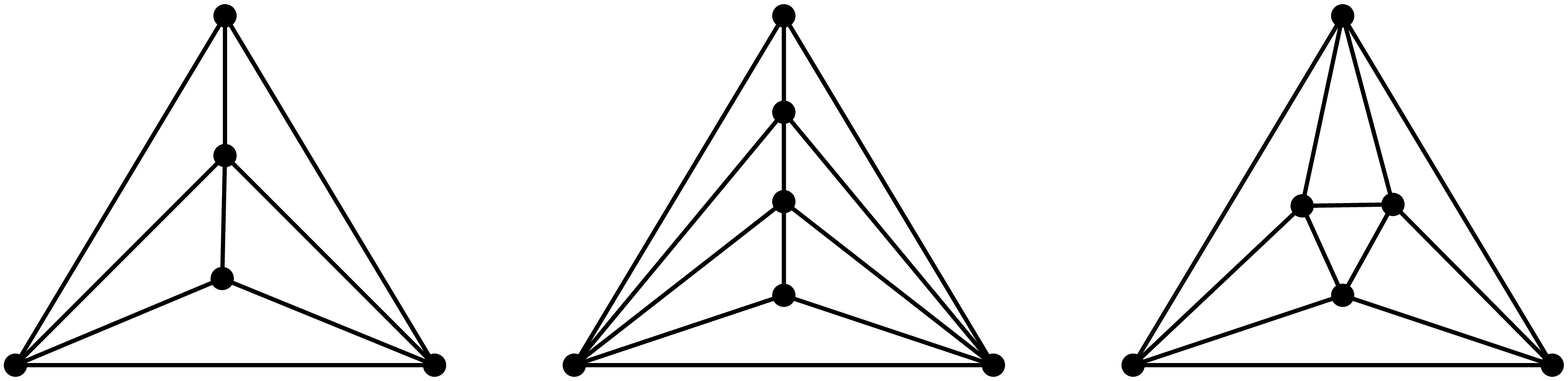}
\caption{All plane triangulations on $5$ and  $6$ vertices.}\label{56}
\end{figure}
We next show that  $ex_{_\mathcal{P}}(n,W_4)=3n-7$ when $  n\in \{5, 6\}$. Note that  all plane triangulations on $n\in\{5,6\}$ vertices are depicted in Figure~\ref{56}, each   containing a copy of $W_4$. Thus, $ex_{_\mathcal{P}}(n,W_4)\le 3n-7$. On the other hand, for all  $  n\in \{5, 6\}$, the planar graph $K_2+(K_2\cup K_{n-4})$ has  $3n-7$ edges  and is   $W_4$-free. Hence,  $ex_{_\mathcal{P}}(n,W_4)= 3n-7$ when   $  n\in \{5, 6\}$. \medskip

\begin{figure}[htbp]
\centering
\includegraphics*[scale=0.35]{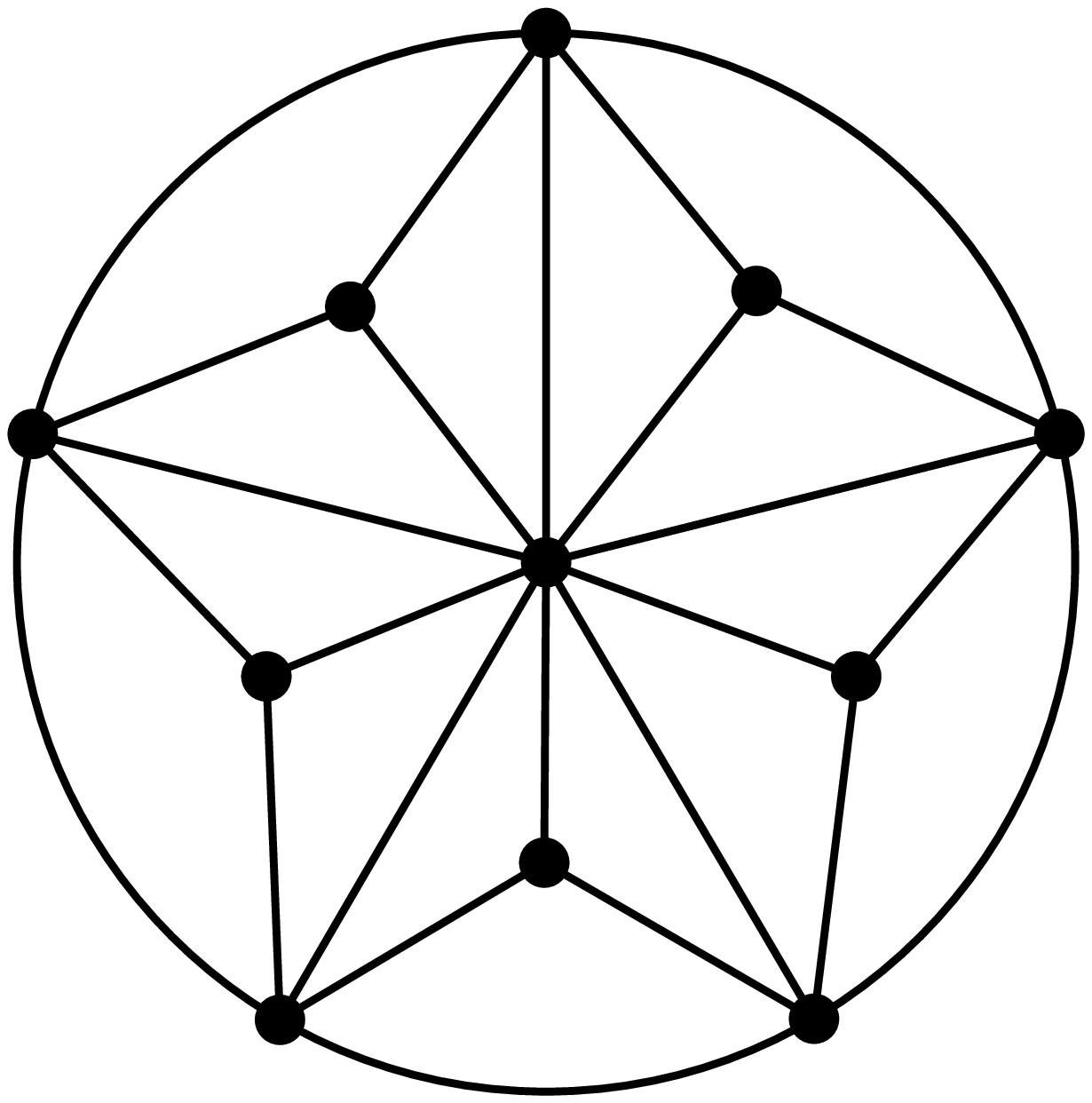}
\caption{Graph  $J$.} \label{T11-2}
\end{figure}

\begin{figure}[htbp]
\centering
\includegraphics*[scale=0.3]{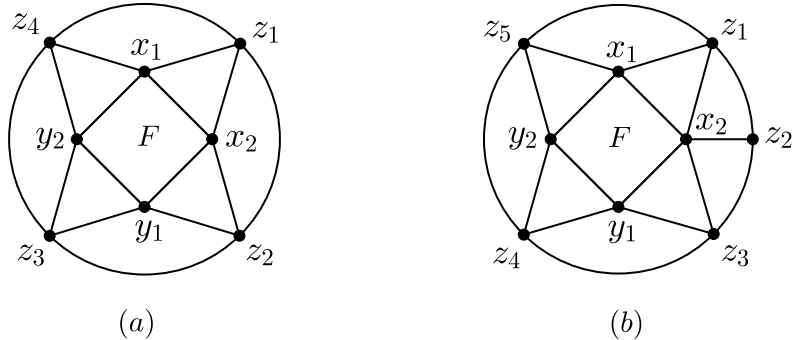}
\caption{The graph $G[A]$.}\label{neighood}
\end{figure}

Finally, we   show that  $ex_{_\mathcal{P}}(n,W_4)=3n-8$  for  all    $n\in\{7,8,9,10,11\}$. The plane graph $J$,  given in Figure~\ref{T11-2},  is $W_4$-free with  $n=11$ vertices and  $3n-8$ edges.
Let $B$ be the set of all vertices of degree $3$ in $J$. Then $|B|=5$. For each $n\in\{7,8,9,10\}$, let $J_n$ be the plane graph obtained from $J$ by deleting $11-n$ vertices in $B$. Then $J_n$ is an induced subgraph of $J$. Clearly,  $J_n$  is $W_4$-free with  $n $ vertices and  $3n-8$ edges.  
Hence, $ex_{_{\mathcal{P}}}(n,W_4)\ge 3n-8$ for all  $n\in\{7,8,9,10,11\}$.
We next show that  $ex_{_{\mathcal{P}}}(n,W_4)\le 3n-8$ for all  $n\in\{7,8,9,10,11\}$.
 Suppose this is not true. Let $G$ be a $W_4$-free  planar graph on $n\in\{7, 8,9,10,11\}$ vertices with  $e(G)\ge 3n-7$.  We choose such a $G$ with $n$ minimum.  Then  $G=T_n$ or  $G=T^-_n$.
Since  each $T_7$,  depicted in  Figure \ref{seven},  contains a copy of $W_4$,  and each $T_7^-$ also contains a copy of $W_4$,  it follows that  $n\in\{8,9,10,11\}$.    Let $u\in V(G)$ with  $d_G(u)=\delta(G)$. Then $\delta(G)\le 4$,  else  $e(G)\ge \frac{5n}{2}>3n-6$ because $n\le11$, a contradiction. Next, if $\delta(G)\le 3$, then $ e(G\less u)\le 3(n-1)-8$ by minimality of $n$ and the fact that $ex_{_{\mathcal{P}}}(n,W_4)\le 3n-8$ when $n=7$. Thus,  $e(G)=e(G\less u)+d_G(u)\le 3(n-1)-8+3=3n-8$, a contradiction. This proves that   $\delta(G)=4$.
Since   $ N_G[u] $ does not  contain    a copy of  $W_{4}$
in $G$, we see that   $G\ne T_n$. Thus   $G=T_n^-$.  We may assume that $G$ is a plane drawing of $T_n^-$  such that the outer face is a $3$-face.   Let $x_1, y_1\in V(G)$ be such that $G+x_1y_1=T_n$.
Then $x_1$ and $y_1$ must lie on the boundary of the unique $4$-face, say $F$,  in $G$.  Let $x_1, x_2, y_1, y_2$ be the vertices on the boundary of $F$ in order. Then $d_G(v)\ge5$ for all $v\in V(G)\less \{x_1, x_2, y_1, y_2\}$,  because $G=T_n^-$ and $ N_G[u] $ does not  contain    a copy of  $W_{4}$
in $G$ for any $u\in V(G)$ with $d_G(u)=4$.   Thus  
 $2(3n-7)=2e(G)\ge 4\cdot 4+5\cdot (n-4)$,
which implies that $n\in\{10,11\}$.  
Suppose  each vertex in $\{x_1, x_2, y_1, y_2\}$ has degree $4$ in $G$. Since  $G=T_n^-$, there must exist four distinct vertices $z_1, z_2, z_3, z_4\in V(G)\less \{x_1, x_2, y_1, y_2\}$ such that $G[A]$ is isomorphic to the graph given in Figure~\ref{neighood}(a),
where $A= \{x_1, x_2, y_1, y_2, z_1, z_2, z_3, z_4\}$. But then  $G$ contains $K_{3,3}$ as a minor,  because $n\in\{10,11\}$ and $d_G(v)\ge5$ for any $v\in V(G)\less A$.
Thus, we may assume that $d_G(x_2)\ge5$. Then $2(3n-7)=2e(G)\ge 4\cdot 3+5\cdot (n-3)$, which implies that $n=11$,  $d_G(v)=4$ for all $v\in\{x_1,   y_1, y_2\}$ and $\Delta(G)=5$.  Thus there exist five distinct vertices $z_1, z_2, z_3, z_4, z_5\in V(G)\less \{x_1, x_2, y_1, y_2\}$ such that $G[A]$ is isomorphic to the graph given in Figure~\ref{neighood}(b), where $A= \{x_1, x_2, y_1, y_2, z_1, z_2, z_3, z_4, z_5\}$. But then $e_{_G}(A, V(G)\less A)=6$, contrary to  $e_{_G}(V(G)\less A, A)\ge 8$ because $n=11$ and $d_G(v)=5$ for any $v\in V(G)\less A$.\medskip

This completes the proof of Theorem~\ref{wheel}.\qed

\section{Proof of Theorem~\ref{K12K2}}   

To establish the desired lower bound, note that the planar graph $K_2+(n-2)K_1$ is $(K_1+2K_2)$-free for all $n\ge5$. Hence,   $ex_{_\mathcal{P}}(n,K_1+2K_2)\ge 2n-3$ for all $n\ge5$. In particular,  $ex_{_\mathcal{P}}(5,K_1+2K_2)\ge 7$.
We next show that  every  $(K_1+2K_2)$-free planar graph on $n\ge 5$ vertices has at most $19n/8-4$ edges.  We proceed the proof by induction on $n$.  Assume first  $n=5$.  Then  $ex_{_\mathcal{P}}(5,K_1+2K_2)= 7 $, because  the only plane triangulation on five vertices,  given in  Figure~\ref{56}(a), is not $(K_1+2K_2)$-free, and any $T_5^-$ is not $(K_1+2K_2)$-free.  Hence, $ex_{_\mathcal{P}}(n,K_1+2K_2)= 7<19n/8-4$ when $n=5$.
So we may assume that $n\ge6$. Let $G$ be a  $(K_1+2K_2)$-free plane graph on $n\ge 6$ vertices. Assume  there exists a vertex $u\in V(G)$ with $d_G(u)\le 2$. By the induction hypothesis,   $e(G\less u)\le 19(n-1)/8-4$ and so $e(G)=e(G\less u)+d_G(u)\le 19(n-1)/8-4+2<19n/8-4$, as desired. So we may assume that $\delta(G)\ge 3$. Next, assume $G$ is disconnected. Then each component of $G$ either is  isomorphic to $K_4$ or has at least  six   
vertices because   $\delta(G)\ge 3$. Let $G_1, \dots, G_p, G_{p+1}, \dots, G_{p+q}$ be all components of $G$ such that $|G_1|=\cdots=|G_p|=4$ and $6\le |G_{p+1}|\le\cdots\le |G_{p+q}|$, where $p\ge 0$ and $q\ge0$ are integers with $p+q\ge2$ and $|G_{p+1}|+\cdots+|G_{p+q}|=n-4p$. Then  $e(G_i)=6$ for all $i\in[p]$,  and  $e(G_j)\le \frac{19|G_j|}{8}-4$ for all $j\in\{p+1, \dots, p+q\}$ by the induction hypothesis. Therefore,
\begin{align*}
e(G)&\le  6p+\frac{19(|G_{p+1}|+\cdots+|G_{p+q}|)}8-4q\\
&=\frac{19n}8-\frac{7p}2-4q \\
&\le \frac{19n}8-\frac{7(p+q)}2 <\frac{19n}8-4,
\end{align*}
as desired.  So we may further assume that $G$ is connected. Then $G$ has no faces of size at most two.
 Hence,
$$2e(G)=3f_3+\sum_{i\ge 4}if_i
  \ge 3f_3+4(f-f_3)
  =4f-f_3,$$
which implies that
\begin{equation}\label{1}f\leq e(G)/2+f_3/4.
\end{equation}
Note that each $3$-vertex is incident with at most three distinct $3$-faces in $G$. Furthermore, since $G$ is   $(K_1+2K_2)$-free, we see that for all $j\ge 4$,  each $j$-vertex is incident with at most two distinct $3$-faces in $G$. Let $U\subseteq V(G)$ denote the set of  $3$-vertices  each  incident with exactly three distinct $3$-faces  in $G$.  Then $U$ must be an independent set in $G$ because  $G$ is connected.  Furthermore, no two vertices in $U$ have a common neighbor in $G$, because $G$ is $(K_1+2K_2)$-free. Thus, $4|U|\le n$ and so $|U|\le n/4$.
 It follows that
\begin{equation}
3f_3\le 3|U|+2(n-|U|)=2n+|U|\le 9n/4, \label{2}
\end{equation}
which implies that    $f_3\le 3n/4$.
This, together with  \eqref{1},  further implies that  $f\le e(G)/2+3n/16$. By Euler's formula, $n-2=e(G)-f\ge e(G)/2-3n/16$. Hence,  $e(G)\le 19n/8-4$, as desired.\medskip

\begin{figure}[htbp]
\centering
\includegraphics*[scale=0.3]{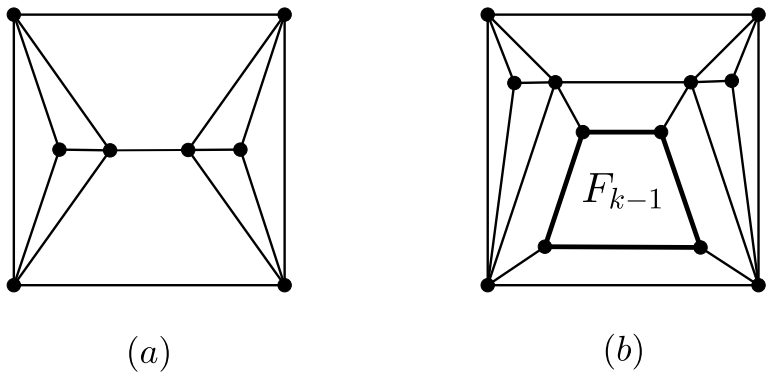}
\caption{The construction of $F_k$.}\label{F}
\end{figure}

From the proof above, we see that  equality in $e(G)\le 19n/8-4$ is achieved for $n$  if and only if equalities  hold  in both  \eqref{1} and \eqref{2},  and in $4|U|\le n$. This implies that $e(G)= 19n/8-4$  for $n$  if and only if $G$ is a connected, $(K_1+2K_2)$-free plane graph on $n$ vertices satisfying: $\delta(G)\ge3$; each $3$-vertex in $G$  is incident with exactly three distinct $3$-faces; each vertex of degree at least $4$ in $G$ is incident with exactly two distinct $3$-faces; each face is either a $3$-face or a $4$-face. We next construct such an extremal plane graph for $n$ and $K_1+2K_2$.  Let $n=8(k+1)$ for some  integer $k\ge0$. Let $F_0$ be the graph depicted in Figure~\ref{F}(a), we then construct $F_k$ of order $n$ recursively  for  all $k\ge1$ via the illustration given in Figure~\ref{F}(b): the entire graph $F_{k-1}$ is placed  into the center  quadrangle of Figure~\ref{F}(b) (in such a way that  the center  bold quadrangle of Figure~\ref{F}(b) is identified with the outer  quadrangle of $F_{k-1}$).
One can check that $F_k$ is $(K_1+2K_2)$-free  with $n=8(k+1)$ vertices and $ 19n/8-4$ edges for all $k\ge0$. \qed\\

\section{Proof of Theorem~\ref{star}} 

\begin{figure}[htbp]
\centering
\includegraphics*[scale=0.22]{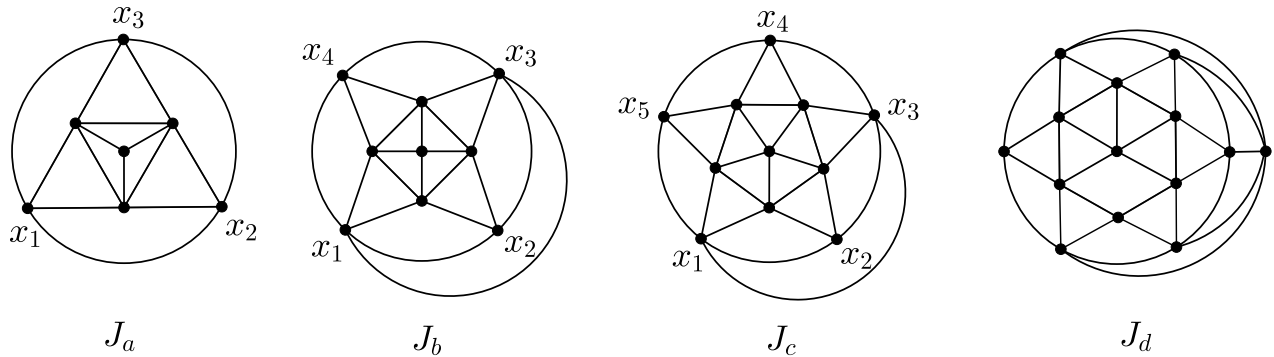}
\caption{The graphs $J_a$, $J_b$, $J_c$ and $J_d$.}\label{7911}
\end{figure}

  By Proposition~\ref{main}(a), $ex_{_\mathcal{P}}(n,K_{1, t})=3n-6$ for all $n\ge t+1\ge 8$. So we may assume that $t\le 6$. We next show that $ex_{_\mathcal{P}}(n,K_{1, 6})=3n-6$ for all  $  n\in\{7,8,9, 10, 12\}$.   Let $J_a,J_b,J_c$  be the plane graphs given in Figure~\ref{7911}. Let  $J'_a$ be the plane triangulation obtained from $J_a$ by adding a new vertex adjacent to $x_1, x_2, x_3$,   $J'_b$  be the plane triangulation obtained from $J_b$ by first deleting the edge $x_1x_3$  and then adding a new vertex adjacent to $x_1, x_2, x_3, x_4$, and $J_c'$ be the plane triangulation obtained from $J_c$ by first deleting the edge  $x_1x_3$  and then adding one new vertex adjacent to $x_1,x_2,x_3,x_4,x_5$. Then  the plane triangulations $J_a$,  $J'_a$,  $J_b$, $J'_b$ and $J_c'$ are $K_{1, 6}$-free   because each of them has maximum degree $5$.
Hence,  $ex_{_\mathcal{P}}(n,K_{1, 6})=3n-6$ for all  $  n\in\{7,8,9, 10, 12\}$.    By Lemma~\ref{Grunbaum1}, no  plane triangulation on $11$ vertices has   maximum degree at most $5$. Hence, every plane triangulation on $n\in\{11, 13, 14\}$ vertices has maximum degree at least $6$. This implies that $ex_{_\mathcal{P}}(n,K_{1, 6})\le3n-7$ for all  $n\in\{11, 13, 14\}$.
Since $J_c$ given in Figure~\ref{7911}  is a $K_{1, 6}$-free  plane graph with $n=11$ vertices and $3n-7$ edges, we have $ex_{_\mathcal{P}}(n,K_{1, 6})=3n-7$ when $n=11$.   By Lemma~\ref{Grunbaum1}, there does not exist any planar  graphs on $n\in\{13,14\}$ vertices with $3n-7$ edges and maximum degree at most $5$.
It follows that  $ex_{_\mathcal{P}}(n,K_{1, 6})\le3n-8$ when $n\in\{ 13, 14\}$. Let   $J_c''$ be the plane graph obtained from $J_c$ by first deleting the edge  $x_1x_3$  and then adding two new adjacent vertices $y_1,y_2$ such that $y_1$ is adjacent to $x_1,x_2,x_3$ and $y_2$ is adjacent to $x_4,x_5$.  Then $J_c''$ and the graph  $J_d$ given in Figure \ref{7911} are $K_{1, 6}$-free plane graph with $n\in\{13,14\}$ vertices and $3n-8$ edges. Hence,  $ex_{_\mathcal{P}}(n,K_{1, 6})=3n-8$ when $n\in\{13,14\}$.  \medskip

 It is easy to see that  $ex_{_\mathcal{P}}(n,K_{1,3})=n$ for all $n\ge 4$, because every $K_{1,3}$-free planar graph has maximum degree at most $2$ and the planar graph $C_n$ is $K_{1,3}$-free with $n$ edges.  We next show that  $ex_{_\mathcal{P}}(n,K_{1,4})=\lfloor3n/2\rfloor$ for all $n\ge 5$. Clearly, $ex_{_\mathcal{P}}(n,K_{1,4})\le \lfloor3n/2\rfloor$ for all $n\ge 5$, because every $K_{1,4}$-free planar graph has maximum degree at most $3$. Next, for all $n\ge 5$, the planar graph obtained from $C_n$ by adding a matching of size $\lfloor n/2\rfloor$ is $K_{1,4}$-free with $\lfloor 3n/2\rfloor$ edges. Hence, $ex_{_\mathcal{P}}(n,K_{1,4})=\lfloor3n/2\rfloor$ for all $n\ge 5$.
 To  determine $ex_{_\mathcal{P}}(n,K_{1, 5})$ for all $n\ge 6$, since  every $K_{1, 5}$-free planar graph on $n\ge 6$ vertices has maximum degree at most $4$, we have $ex_{_\mathcal{P}}(n,K_{1, 5})\le 2n $ for all $n\ge 6$. Let $J_a''$ be the plane triangulation obtained from $J_a$ by deleting the unique $3$-vertex. Since $J_a''$ is $K_{1,5}$-free plane graph on $n=6$ vertices with $2n$ edges, we have $ex_{_\mathcal{P}}(n,K_{1, 5})=2n $ when $n=6$. By Lemma~\ref{Grunbaum1}, no  planar graph with  $n=7$ vertices and  $2n$ edges has  maximum degree at most $4$. Hence, $ex_{_\mathcal{P}}(n,K_{1, 5})\le 2n-1 $ when $n=7$. Let $J_a'''$ be the plane graph obtained from $J_a''$ by first deleting  the edge $x_1x_2$ and then adding a new vertex adjacent to $x_1, x_2$ only. Note that $J_a'''$ is a $K_{1,5}$-free plane graph on $n=7$ vertices with $2n-1$ edges, we see that $ex_{_\mathcal{P}}(n,K_{1, 5})=2n-1$ when $n=7$. Next, for all $n\ge 8$, let $C$ be a cycle on $2\lfloor n/2\rfloor$ vertices with vertices $u_1,\ldots,u_{\lfloor \frac{n}{2}\rfloor}, w_{\lfloor \frac{n}{2}\rfloor},\ldots,w_1$ in order. Let $H$ be  the plane graph obtained from $C$ by adding the path with vertices $w_1,u_2,w_2,u_3,\ldots,w_{\lfloor \frac{n}{2}\rfloor-1},u_{\lfloor \frac{n}{2}\rfloor}$ in order. When $n$ is even, the planar graph $H+u_1u_{\lfloor \frac{n}{2}\rfloor}+u_1w_{\lfloor \frac{n}{2}\rfloor}+w_1w_{\lfloor \frac{n}{2}\rfloor}$ is $K_{1,5}$-free with $2n$ edges. When $n$ is odd, let $H'$ be obtained from $H$ by first deleting the edge $u_2u_3$ and then adding a new vertex $u$ adjacent to $u_2$ and $u_3$. Then the planar graph $H'+uu_1+uu_{\lfloor \frac{n}{2}\rfloor}+w_1w_{\lfloor \frac{n}{2}\rfloor}+u_1w_{\lfloor \frac{n}{2}\rfloor}$ is $K_{1,5}$-free with $2n$ edges. Hence, $ex_{_\mathcal{P}}(n,K_{1, 5})=2n $ for all $n\ge 8$.\medskip

    It remains to show  that $ex_{_\mathcal{P}}(n,K_{1, 6})= \lfloor5n/2\rfloor $     for all    $n\ge 15$.
Clearly, $ex_{_\mathcal{P}}(n,K_{1, 6})\le  \lfloor5n/2\rfloor $  for all    $n\ge 15$, because every $K_{1, 6}$-free planar graph on $n\ge 15$ vertices has maximum degree at most $5$. We next show that $ex_{_\mathcal{P}}(n,K_{1, 6})\ge  \lfloor5n/2\rfloor $  for all    $n\ge 15$.  Let $n: =4q+r\ge15$, where $q\ge3$ and $r\in\{0,1,2,3\}$.  Let $p\in\{q, q+1\}$. Let $C^1 $ and  $C^2$  be two vertex-disjoint cycles with vertices $x_1, x_2, \ldots x_q$ in order and  $y_1, y_2, \ldots y_{p}$ in order, respectively. Let $C^3$ be a cycle of length $q+p$  with vertices $b_1, a_1, b_2, a_2, \ldots b_q, a_q$ in order when $p=q$,  and  $b_1, a_1, b_2, a_2, \ldots b_q, a_q, b_{q+1}$ in order when $p=q+1$.
Let $R_p$  be the plane graph on $2q+2p$ vertices obtained from disjoint copies of $C^1$, $C^2$ and $C^3$ by making $x_i$ adjacent to  $\{a_i,b_i,b_{i+1}\}$, and $y_j$   adjacent to $\{b_j,a_{j-1},a_j\}$ for all $1\le i\le q$ and $1\le j\le p$,
\begin{figure}[htbp]
\centering
\includegraphics*[scale=0.28]{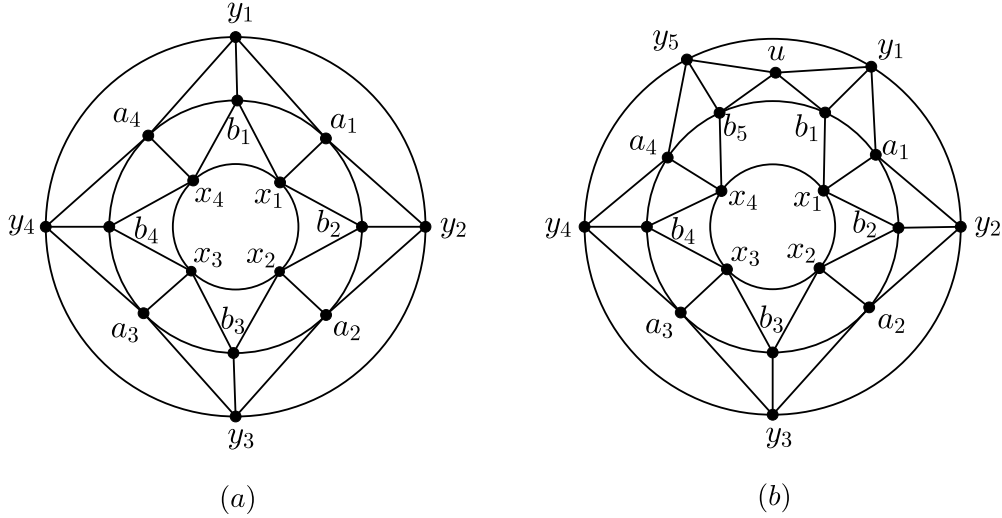}
\caption{Almost $5$-regular plane graphs on $4q+r$ vertices when $q=4$ and $r\in\{0,3\}$.}\label{R4}
\end{figure}
 where  all arithmetic on the indices $i+1$ and $j-1$ here are done modulo $p$. Then $R_p$ is $K_{1,6}$-free planar graph with $10q$ edges when $p=q$ and  $10q+3$ edges when $p=q+1$. The construction of   $R_p$ on $4q$ vertices when $q=4$ is depicted in Figure~\ref{R4}(a). When $n=4q$ for $q\ge 4$, the planar graph $R_p$ with $p=q$ is $K_{1,6}$-free with $10q$ edges and so $ex_{_\mathcal{P}}(n,K_{1,6})=10q=5n/2$.   When $n=4q+1$ for $q\ge 4$, let $R^1$ be  obtained from $R_p-y_2y_3-y_1y_p$ with $p=q$ by adding a new vertex $u$ adjacent to $y_2$ and $y_3$. Then the planar graph $R^1+uy_1+uy_p$ is $K_{1,6}$-free with $10q+2$ edges. Hence, $ex_{_\mathcal{P}}(n,K_{1,6})=10q+2=\lfloor 5n/2\rfloor$ when $n=4q+1$ for $q\ge 4$. When $n=4q+2$ for $q\ge 4$, let $R^2$ be  obtained from $R_p-y_2y_3-y_1y_p-x_2x_3-x_1x_p-b_1a_p$ with $p=q$ by adding two new vertices $u$ adjacent to $y_2$ and $y_3$ and $v$ adjacent to $x_2$ and $x_3$. Then the planar graph $R^2+uy_1+uy_p+ua_q+vx_1+vx_q+vb_1$ is $K_{1,6}$-free with $10q+5$ edges. Hence, $ex_{_\mathcal{P}}(n,K_{1,6})=10q+5=5n/2$ when $n=4q+2$ for $q\ge 4$. When $n=4q+3$ for $n\ge 3$, the planar graph obtained from $R_p$  with $p=q+1$ by adding a new vertex $u$ adjacent to $y_1,b_1,y_p,b_p$, given in Figure~\ref{R4}(b) when $q=4$,  is $K_{1,6}$-free with $10q+7$ edges. Hence, $ex_{_\mathcal{P}}(n,K_{1,6})=10q+7=\lfloor5n/2\rfloor$ when $n=4q+3$ for $q\ge 3$.

\bigskip

\section{Proof of Theorem~\ref{K13K2}} 

 Since the plane triangulations $J_a$,   $J'_a$, $J_b$ and $J'_b$ constructed in the proof of Theorem~\ref{star} is  $(K_1+3K_2)$-free,
we see that  $ex_{_\mathcal{P}}(n,K_1+3K_2)=3n-6$ for all  $n\in\{7,8,9,10\}$. To determine $ex_{_\mathcal{P}}(11,K_1+3K_2)$, note that the plane graph $J_c$ given in Figure~\ref{7911} with  $n=11$ vertices and $3n-7$ edges  is $(K_1+3K_2)$-free. Thus $ex_{_\mathcal{P}}(n,K_1+3K_2)\ge 3n-7$ when $n=11$.   By Lemma~\ref{Grunbaum1}, no plane triangulation on $11$ vertices has  maximum degree at most $5$. Hence, every plane triangulation on $11$ vertices must contain a vertex of degree at least $6$ (and so contains a copy of $K_1+3K_2$),  which implies that $ex_{_\mathcal{P}}(n,K_1+3K_2)= 3n-7$ when $n=11$. Since every $K_{1,6}$-free graph is certainly $(K_1+3K_2)$-free, 
 by Theorem \ref{star}, $ex_{_\mathcal{P}}(n,K_1+3K_2)= ex_{_\mathcal{P}}(n,K_{1,6}) =3n-6$ when $n= 12$,   $ex_{_\mathcal{P}}(n,K_1+3K_2)\ge ex_{_\mathcal{P}}(n,K_{1,6})=3n-8$ when $n\in\{ 13, 14\}$, and $ex_{_\mathcal{P}}(n,K_1+3K_2)\ge ex_{_\mathcal{P}}(n,K_{1,6})=\lfloor5n/2\rfloor $ when $n\ge 15$. Since every plane triangulation on $n\in\{13,14\}$ vertices has maximum degree at least $6$, we see that $ex_{_\mathcal{P}}(n,K_1+3K_2)\le 3n-7$ when $n\in\{13,14\}$. By Lemma~\ref{Grunbaum1}, every  $T_n^-$ with $n\in\{13,14\}$   has maximum degree at least $6$ and so contains a copy of $K_1+3K_2$. It follows that $ex_{_\mathcal{P}}(n,K_1+3K_2)=3n-8<  17n/6 -4$ when $n\in\{13,14\}$. \medskip

    We next show that  every  $(K_1+3K_2)$-free planar graph $G$ on $n\ge 13$ vertices has at most $17n/6-4$ edges.
We proceed the proof by induction on $n$.  This is trivially true when  $n\in\{13,14\}$. So we may assume that $n\ge15$.    Assume  there exists a vertex $u\in V(G)$ with $d_G(u)\le 2$. By the induction hypothesis, $e(G\less u)\le 17(n-1)/6 -4$ and so $e(G)=e(G\less u)+d_G(u)\le 17(n-1)/6 -4+2<17 n/6-4$, as desired. So we may assume that $\delta(G)\ge 3$.
Assume next that  $G$ is disconnected.
Let $G_1, \dots, G_p, G_{p+1}, \dots, G_{p+q}$ be all components of $G$ such that $|G_1|\le\cdots\le|G_p|\le 12$ and $13\le |G_{p+1}|\le\cdots\le |G_{p+q}|$, where $p\ge 0$ and $q\ge0$ are integers with $p+q\ge2$ and $|G_{1}|+\cdots+|G_{p+q}|=n$. Then  $e(G_i)\le3|G_i|-6$ for all $i\in[p]$,  and  $e(G_j)\le  17|G_j|/6-4$ for all $j\in\{p+1, \dots, p+q\}$ by the induction hypothesis. Therefore,
\begin{align*}
e(G)&\le 3(|G_{1}|+\cdots+|G_{p}|)-6p+\frac{17(|G_{p+1}|+\cdots+|G_{p+q}|)}6-4q\\
&=\frac{17n}6-(6p+4q)+\frac{|G_{1}|+\cdots+|G_{p}|}6\\
&\le \frac{17n}6-(6p+4q)+2p=\frac{17n}6-4(p+q)<\frac{17n}6-4,
\end{align*}
as desired.  So we may further assume that $G$ is connected. Then $G$ has no faces of size at most two. Hence,
$$2e(G)=3f_3+\sum_{i\geq 4}if_i
\geq 3f_3+4(f-f_3)
=4f-f_3,$$
which implies that $f\leq e(G)/2+f_3/4$.
Note that $  n_{_3}(G)\ge0$ and $n_{_5}(G)< n$; and  for all $i\in \{3,4, 5\}$, each $i$-vertex
is incident with  at most $i$ $3$-faces. Furthermore,     for all $j\ge6$, each $j$-vertex   is incident with  at most four $3$-faces  because $G$ is $(K_1+3K_2)$-free  and $n\ge 15$. It follows that
\begin{align*}
3f_3\le 3n_{_3}(G)+4n_{_4}(G)+5n_{_5}(G)+4(n-n_{_3}(G)-n_{_4}(G)-n_{_5}(G))
=4n-n_{_3}(G)+n_{_5}(G)< 5n,
\end{align*}
which implies that $f_3< 5n/3$. This, together with    the fact that $f\leq e(G)/2+f_3/4$,  further implies that   $f< e(G)/2+5n/12$. By Euler's formula, $n-2=e(G)-f> e(G)/2-5n/12$. Hence,    $e(G)<17n/6-4$. \qed

\section{Concluding  remarks}

\begin{figure}[htbp]
\centering
\includegraphics*[scale=0.28]{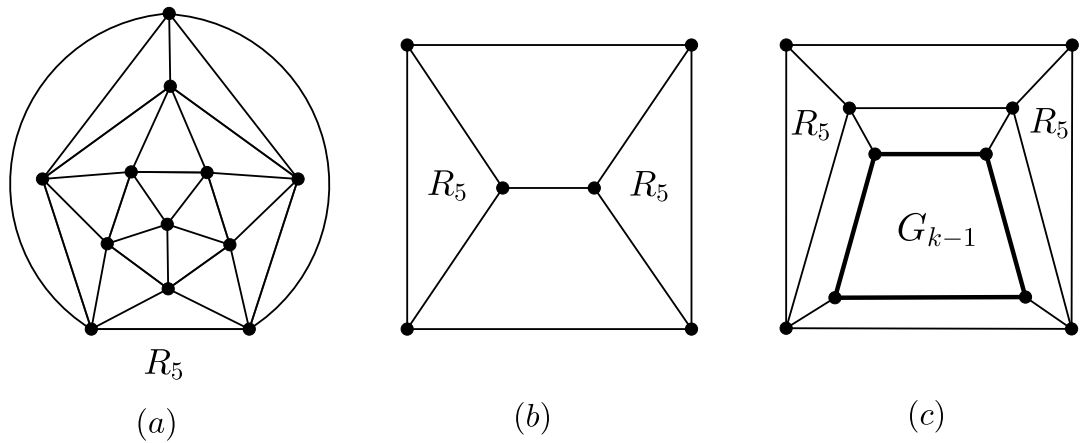}
\caption{The construction of $G_k$.}\label{lower}
\end{figure}

The lower bound in Theorem~\ref{K13K2} can be further improved when $n$ is divisible by $24$. To see this,  let $n=24(k+1)$ for some  integer $k\ge0$.   Let $R_5$ be the $5$-regular plane graph on twelve  vertices  given in Figure~\ref{lower}(a),  and let $G_0$ be the plane graph obtained from two  disjoint copies of $R_5$ by adding three independent edges between their outer faces,  as depicted  in Figure~\ref{lower}(b).  We construct $G_k$ of order $n$ recursively  for  all $k\ge1$ via the illustration given in Figure~\ref{lower}(c): the entire graph $G_{k-1}$ is placed  into the center  quadrangle of Figure~\ref{lower}(c) (in such a way that  the center  bold  quadrangle of Figure~\ref{lower}(c) is identified with the outer  quadrangle of  $G_{k-1}$).
One can check that $G_k$ is $(K_1+3K_2)$-free  with $n=24(k+1)$ vertices and $ 67n/24-4$ edges for all $k\ge0$.

\bigskip
As mentioned earlier, it seems non-trivial  to determine $ex_{_\mathcal{P}}(n,H)$ for all     $K_4$-free planar graphs $H$  with exactly one vertex, say $u$, satisfying $d_H(u)=\Delta(H)\le6$ and $\Delta(H[N(u)])\le2$.   We conclude this section by giving an upper bound (but not tight) for $ex_{_\mathcal{P}}(n,K_1+H)$, where    $H$ is  a  disjoint union of paths. 

\begin{thm}\label{K1Pt}
Let $4\le t\le6$ be an integer and let $H$ be a graph on $t$ vertices such that $H$ is a  disjoint union of paths. Then $ex_{_\mathcal{P}}(n,K_1+H)\le  \frac{13(t-1)n}{4t-2}-\frac{12(t-1)}{2t-1}$ for all $n\ge t+1$.
\end{thm}

\pf Let $t$ and $H$ be given as in the statement. Since $H$ is a subgraph of $K_1+P_t$, it suffices to show that $ex_{_\mathcal{P}}(n,K_1+P_t)\le \frac{13(t-1)n}{4t-2}-\frac{12(t-1)}{2t-1}$ for all $n\ge t+1$.     Let $G$ be a $(K_1+P_{t})$-free planar graph on $n\ge t+1$ vertices.  We next show that $e(G)\le  \frac{13(t-1)n}{4t-2}-\frac{12(t-1)}{2t-1}$ by induction on $n $. This is trivially true when $n=t+1$ because $e(G)\le 3(t+1)-6\le  \frac{13(t-1)(t+1)}{4t-2}-\frac{12(t-1)}{2t-1}$  for $5\le t\le 6$ and $e(G)\le 3(t+1)-7\le  \frac{13(t-1)(t+1)}{4t-2}-\frac{12(t-1)}{2t-1}$ for $t=4$ . So we may assume that $n\ge t+2$.  We may further  assume that $\delta(G)\ge 3$ and $G$ is connected. Hence, $$2e(G)=3f_3+\sum_{i\geq 4}if_i
\geq 3f_3+4(f-f_3)
=4f-f_3,$$
which implies that $f\leq e(G)/2+f_3/4$.
Note that for all $3\le i\le t-1$, each $i$-vertex
is incident with  at most $i$ many $3$-faces, and for all $j\ge t$, each $j$-vertex   is incident with  at most $(t-2)\lfloor\frac{j}{t-1}\rfloor$   many $3$-faces,   because $G$ is $(K_1+P_t)$-free. It follows that
\begin{align*}
3f_3&\le \sum_{i=3}^{t-1}i\cdot n_i(G)+\sum_{j\ge t}(t-2)\left\lfloor\frac{j}{t-1}\right\rfloor \cdot n_j(G)\\
&\le\sum_{i=3}^{t-1}\left(\frac{i\cdot n_i(G)}{t-1}+\frac{(t-2)i\cdot n_i(G)}{t-1}\right)+\frac{t-2}{t-1}\sum_{j\ge t}j\cdot n_j(G)\\
&=\sum_{i=3}^{t-1}\frac{i\cdot n_i(G)}{t-1} +\frac{t-2}{t-1}\sum_{\ell\ge 3}\ell\cdot n_\ell(G)\\
&= \sum_{i=3}^{t-1}\frac{i\cdot n_i(G)}{t-1}+\left(\frac{t-2}{t-1}\cdot2e(G)\right)\\
&=  \sum_{i=3}^{t-1}n_i(G)-\sum_{i=3}^{t-1}\frac{(t-1-i)\cdot n_i(G)}{t-1}+\left(\frac{t-2}{t-1}\cdot2e(G)\right)\\
&<n+\frac{t-2}{t-1}\cdot2e(G), \label{a}
\end{align*}
which implies that $f_3< \frac{n}{3}+\frac{2(t-2)}{3(t-1)}\cdot e(G)$. This, together with    the fact that $f\leq e(G)/2+f_3/4$,  further implies that   $f\le \frac{e(G)}{2}+\frac{n}{12}+\frac{(t-2)}{6(t-1)}\cdot e(G)=\frac{(4t-5)}{6(t-1)}\cdot e(G)+\frac{n}{12}$. By Euler's formula, $n-2=e(G)-f\ge \frac{(2t-1)}{6(t-1)}\cdot e(G)-\frac{n}{12}$. Hence,    $e(G)\le \frac{13(t-1)n}{4t-2}-\frac{12(t-1)}{2t-1}$. \qed\\

\noindent {\bf Acknowledgments.}
\noindent  Zi-Xia Song would like to thank Yongtang Shi and  the Chern Institute of Mathematics at Nankai University for hospitality and support during her visit  in June 2018. \medskip

Yongxin Lan and Yongtang Shi are partially supported by the National Natural Science Foundation
of China and the Natural Science Foundation of Tianjin (No.17JCQNJC00300).


\frenchspacing
\begin{thebibliography}{14}
\bibitem{2013Bollobas}
B. Bollob\'{a}s,   Modern Graph Theory, Springer, 2013.
%
\vspace {-0.25cm}
%
\bibitem{CGPW2003}
G. Chen, R. J. Gould, F. Pfender and B. Wei, Extremal graphs for intersecting cliques,  J. Combin. Theory Ser. B 89 (2003), 159--171.

\vspace{-0.25cm}


\bibitem{Dowden2016}
C. Dowden, Extremal $C_4$-free/$C_5$-free planar graphs,  J. Graph Theory  83 (2016),  213--230.
%
\vspace {-0.25cm}

\bibitem{Dzido2013}
T. Dzido, A note on Tur\'an numbers for even wheels, Graphs Combin. 29 (2013), 1305--1309.

\vspace{-0.25cm}

\bibitem{DJ2017}
T. Dzido and A. Jastrz\c{e}bski, Tur\'an numbers for odd wheels,  Discrete Math. 341 (2018) 1150--1154.

\vspace{-0.25cm}


\bibitem{EFGG1995}
P. Erd\H{o}s, Z. F\"{u}redi, R. J. Gould and D. S. Gunderson, Extremal graphs for intersecting triangles,   J. Combin. Theory Ser. B 64 (1995), 89--100.

\vspace{-0.25cm}

\bibitem{ErdosStone} P.  Erd\H{o}s and A. Stone, On the structure of linear graphs, Bull. Amer. Math. Soc. 52 (1946), 1087--1091.
 \vspace {-0.25cm}
\bibitem{Furedi} Z. F\"uredi,  Tur\'an type problems, ``Surveys in Combinatorics", London Math. Soc. Lecture Note Ser. 166, Cambridge Univ. Press, Cambridge, 1991, pp. 253--300.
%
\vspace {-0.25cm}
%
\bibitem{FurediJiang} Z. F\"uredi and T. Jiang, Hypergraph Tur\'an numbers of linear cycles, J. Combin. Theory Ser. A 123 (2014),  252--270.
%
\vspace {-0.25cm}
%
\bibitem{FurediJiangSeiver} Z. F\"uredi, T. Jiang and R. Seiver, Exact solution of the hypergraph Tur\'an problem for $k$-uniform linear paths, Combinatorica 34 (2014),  299--322.
%
\vspace {-0.25cm}
%
\bibitem{Keevash} P.  Keevash,  Hypergraph Tur\'an problems,   ``Surveys in Combinatorics", London Math. Soc. Lecture Note Ser. 392, Cambridge Univ. Press, Cambridge, 2011, pp. 83--139.
%
\vspace {-0.25cm}

\bibitem{Kostochka1} A. Kostochka,  D.  Mubayi and J. Verstra\"ete, Tur\'an problems and shadows I: Paths and cycles, J. Combin. Theory Ser. A 129 (2015), 57--79.
%
\vspace {-0.25cm}

\bibitem{LSS}
Y. Lan, Y. Shi and Z-X. Song, Planar Tur\'an numbers for Theta graphs and paths of small order, arXiv:1711.01614, 2017.
\vspace {-0.25cm}
 
 

\bibitem{SH1977}
E. F. Schmeichel and S. L. Hakimi, On planar graphical degree sequences,   SIAM J. Appl. Math. 32 (1977), 598--609.

\vspace {-0.25cm}

\bibitem{Turan1941}
P.  Tur\'an, On an extremal problem in graph theory,  Mat. Fiz. Lapok. 48 (1941), 436--452.


 

\end{thebibliography}
\end{document}